\documentclass[12pt,fullpage,report]{amsart}
\usepackage{amsfonts}
\usepackage{amsmath}
\usepackage{amssymb}

\newcommand{\reals}{{\mathbb{R}}}
\newcommand{\vareps}{{\varepsilon}}
\newcommand{\eps}{{\varepsilon}}
\newcommand{\Lp}[2]{\|{#1}\|_{L^{#2}}}
\newcommand{\kovrtwo}{{\left\lceil\frac{k}{2}\right\rceil}}
\newcommand{\half}{{\frac{1}{2}}}
\newcommand{\scriptj}{{\mathcal{J}}}
\newcommand{\scriptt}{{\mathcal{T}}}
\newcommand{\scripts}{{\mathcal{S}}}
\newcommand{\scriptk}{{\mathcal{K}}}
\newcommand{\scriptr}{{\mathcal{R}}}
\newcommand{\scriptb}{{\mathcal{B}}}

\newtheorem{theorem}{Theorem}[section]
\newtheorem{lemma}[theorem]{Lemma}
\newtheorem*{hyp1*}{Hypothesis 1}
\newtheorem*{hyp2*}{Hypothesis 2}

\title{Endpoint bounds for a generalized Radon transform}
\thanks{The final version of this article will appear in the J.\ of the London Math.\ Soc.}

\author{Betsy Stovall}
\address{Department of Mathematics, UC Berkeley, Berkeley, CA 94720-3840}
\email{betsy@@math.berkeley.edu}
\thanks{The author was supported in part by NSF grant DMS-040126.}


\begin{document}
\maketitle

\begin{abstract}
We prove that convolution with affine arclength measure on the curve parametrized by $h(t) := (t,t^2,\ldots,t^n)$ is a bounded operator from $L^p(\reals^n)$ to $L^q(\reals^n)$ for the full conjectured range of exponents, improving on a result due to M.~Christ.  We also obtain nearly sharp Lorentz space bounds.
\end{abstract}

\section{Introduction}

Let $T$ be the operator defined on Borel measurable functions on
$\reals^n$ by $Tf(x) = \int_{-1}^1 f(x-h(t))dt$, where $h(t) =
(t,t^2,\ldots,t^n)$. The study of $L^p \to L^q$ bounds for this
operator was initiated by Littman in \cite{Litt}; there it was
proved that when $n=2$, $T$ extends as a bounded operator from $L^p$
to $L^q$ if and only if $(p^{-1},q^{-1})$ lies in the convex hull of
the points $(0,0),(1,1),(2/3,1/3)$.  Later, Oberlin addressed the
case $n=3$ in \cite{Obe1}, proving that $T$ is bounded from $L^p$ to
$L^q$ if and only if $(p^{-1},q^{-1})$ belongs to the convex hull of
the points $(0,0),(1,1),(1/2,1/3),(2/3,1/2)$.

In general, we let $p_n=\frac{n+1}{2}$ and
$q_n=\frac{n+1}{2}\frac{n}{n-1}$ and let $\scriptr_n$ be the convex
hull of the points $(0,0),(1,1),(p_n^{-1},q_n^{-1}),(1-q_n^{-1},
1-p_n^{-1})$.  In \cite{CCC}, Christ used combinatorial methods to
prove that when $n \geq 4$, $T$ is of restricted weak type
$(p_n,q_n)$, which by interpolation and duality proved that $T$ maps
$L^p(\reals^n)$ to $L^q(\reals^n)$ if $(p^{-1},q^{-1})$ lies in
$\scriptr_n \backslash \{(p_n^{-1},q_n^{-1}),(1-q_n^{-1},
1-p_n^{-1})\}$.

Using techniques developed by Christ in \cite{CCC} and \cite{QEx},
we prove that when $n \geq 2$, $T$ maps $L^{p_n}(\reals^n)$ to
$L^{q_n}(\reals^n)$; we also obtain an improvement in Lorentz
spaces.

\begin{theorem}
For $n \geq 2$, $T$ extends as a bounded operator from
$L^p(\reals^n)$ to $L^q(\reals^n)$ whenever $(p^{-1},q^{-1})$ lies
in $\scriptr_n$.  Moreover, $T$ maps $L^{p_n,u}(\reals^n)$ boundedly
into $L^{q_n,v}(\reals^n)$ and $L^{q_n',v'}(\reals^n)$ into
$L^{p_n',u'}(\reals^n)$ whenever $u<q_n$, $v>p_n$, and $u<v$.
\end{theorem}

As mentioned above, when $n=2,3$, boundedness at the Lebesgue
endpoints is already known.  The Lorentz space bounds attained here
are known when $n=2$.  These bounds were first shown in \cite{BOS} (including endpoints).  In \cite{QEx} there is an alternative argument for the Lorentz bounds, along the lines of those here (indeed, the $n=2$ case here and there are identical).  In the introduction of \cite{QEx}, Christ outlines an argument which produces Lorentz space bounds when the Lebesgue space exponents are integers. In a recent paper
\cite{BandS}, Bennett and Seeger have shown that when $n=3$, $T$
maps $L^{3/2,2}$ boundedly into $L^{2}$ (and hence $L^{2} \to L^{3,2}$).  More recently, in \cite{DLW} Dendrinos, Laghi, and Wright have established the analogue of our theorem for convolution with affine arclength measure along arbitrary polynomial curves in low dimensions.

We do not address the Sobolev regularity  of this operator.  See \cite{PS}, for some recent results in that direction.

In \S \ref{sec:sharp}, we will show that if $T$ is bounded from
$L^{p_n,u}$ to $L^{q_n,v}$, then the inequalities $u \leq v$, $u
\leq q_n$, and $v \geq p_n$ must hold, so this result is sharp up to
Lorentz space endpoints.  The $L^2(\reals^3) \to
L^{3/2,2}(\reals^3)$ bound obtained by Bennett and Seeger indicates
that this result is still not optimal, but the author has not been able to
extend this proof to the Lorentz space endpoints.

Some work on related operators has been carried out by Tao and Wright in
\cite{TW}, Christ in \cite{ChRL}, and Gressman in \cite{Gress_poly}, for instance. In
\cite{TW}, Tao and Wright considered a far more general class of
operators defined by integration along smoothly varying families of
curves, proving Lebesgue space bounds which are sharp up to
endpoints.  Using partially alternative techniques, Christ reproved the same
bounds in \cite{ChRL}.  Since the methods used here rely heavily on
the polynomial structure of the operator $T$ they do not seem to
generalize to the $C^{\infty}$ case considered by those authors.  In the polynomial case of the Tao-Wright theorem, the restricted weak-type bounds at the endpoints have been proved by Gressman in \cite{Gress_poly}.

This author hopes that with some modifications, the argument in this paper can be used to establish strong-type endpoint bounds (and an improvement in Lorentz spaces) for a more general class of polynomial curves, such as that in \cite{Gress}.

As mentioned above, a quite recent result of Dendrinos, Laghi, and Wright in \cite{DLW} (the authors also use the methods of \cite{QEx}) establishes sharp Lebesgue space bounds (with an accompanying Lorentz space improvement) for convolution with affine arclength measure along polynomial curves in dimensions 2 and 3.  

\subsection*{Acknowledgements}  
The author would like to thank her advisor, Michael Christ, for suggesting this problem and for his advice and help throughout this project.  The author would also like to thank the anonymous reviewer at the JLMS for comments which she believes led to significant improvements in the exposition.

\section{On notation and other preliminary remarks}
\label{sec:prelim}

\subsection*{Notation} Most of the notation we will use is fairly
standard.  If $1 \leq p \leq \infty$, we denote by
$p'$ the exponent dual to $p$.  We use $|\cdot|$ to indicate Lebesgue measure and $\#$ for cardinality.  When $A$ and $B$ are
non-negative real numbers, we write $A \lesssim B$ to mean $A \leq
CB$ for an implicit constant $C$, and $A \sim B$ when $A \lesssim B$ and $B \lesssim A$.  In addition, for $x$ a real number,
$\lceil x \rceil$ and $\lfloor x \rfloor$ are the least integer
greater than or equal to and the greatest integer less than or equal
to x, respectively.  We will also employ the somewhat less standard
notation $\scriptt(E,F):=\langle T
\chi_E,\chi_F\rangle$ when $E$ and $F$ are Borel sets and $T$ is a
linear operator.

\subsection*{The endpoint $(q_n',p_n')$} For the remainder of the paper we
will focus on $L^{p_n,u} \to L^{q_n,v}$ bounds (and
counter-examples), as these imply $L^{q_n',v'} \to L^{p_n',u'}$
bounds (and counter-examples) by duality and the fact that $T^*$ is
essentially the same operator as $T$.

\subsection*{A related operator} We note here that if $T$ is bounded from
$L^{p_n,u}$ to $L^{q_n,v}$ and $0<R<\infty$ then the operator $T_R$
defined by
\[
    T_Rf(x) = \int_{-R}^{R}f(x-h(t))dt
\]
is also bounded from $L^{p_n,u}$ to $L^{q_n,v}$, with a bound
independent of $R$. To see this, note first that $D_R^{-1} \circ T
\circ D_R = R^{-1}T_R$, where $D_R$ is the anisotropic scaling of
$\reals^n$ defined by
\[
    D_R(x_1,x_2,\ldots,x_n) = (Rx_1,R^2x_2,\ldots,R^nx_n),
\]
and second that any $L^{p_n,u} \to L^{q_n,v}$ bound scales under
this transformation.  From this, the operator $T_{\infty}$ is also
bounded from $L^{p_n,u}$ to $L^{q_n,v}$.  By duality and
interpolation, $T_{\infty}$ maps $L^p$ boundedly into $L^q$ whenever
$(p^{-1},q^{-1})$ lies in the line segment $[(p_n^{-1},q_n^{-1}),
((q_n')^{-1},(p_n')^{-1})]$.

\subsection*{Outline}  In \S\ref{sec:sharp}, we will show that our result is
sharp up to endpoints and review the argument that
$T:L^p \to L^q$ is bounded only if $(p^{-1},q^{-1}) \in \scriptr_n$.
In \S\ref{sec:reduction}, we leave the setting of our particular
operator $T$ and state two hypotheses--essentially multilinear bounds
involving characteristic functions of sets--which suffice to prove
$L^{r,u} \to L^{s,v}$ bounds for the operator $T$.  A proof of this fact, using an argument developed in \cite{QEx}, will be postponed until the appendix.  Finally, in
\S\ref{sec:multilinear}, we prove that the hypotheses from
\S\ref{sec:reduction} do in fact hold.  For this, we use an
iteration scheme and ``band structure" argument similar to that in
\cite{CCC}.

\section{Almost sharpness} \label{sec:sharp}

Given $\vareps>0$, we let $N_{\vareps}$ be the
$\vareps$-neighborhood of the curve $-h([-1,1])$.  For $r>0$, we let
$N_{\vareps,r}=D_r(N_{\vareps})$, where $D_r$ is the anisotropic scaling from 
\S\ref{sec:prelim}. We also define $B_{\vareps}$ to be the
$\vareps$-neighborhood of 0 and let $B_{\vareps,r} =
D_r(B_{\vareps})$.

To see that $T$ can only be of restricted weak-type $(p,q)$ when
$(p^{-1},q^{-1}) \in \scriptr_n$, one need only compare
$\scriptt(E,F)$ and $|E|^{1/p}|F|^{1/q'}$ for the pairs
$E=N_{\vareps,r}, F=B_{\vareps,r}$ and $E=B_{\vareps,r},
F=-N_{\vareps,r}$ when $0<\vareps,r<1$, and $E=B_R,F=B_R$, when
$R>1$.  See \cite{CCC}.

If $x \in \reals^n$, we define the translates
$N_{\vareps,r}(x):=N_{\vareps,r}+\{x\}$ and $B_{\vareps,r}(x):=
B_{\vareps,r}+\{x\}$.  We will show that if $T$ is a bounded
operator between Lorentz spaces $L^{p_n,u}$ and $L^{q_n,v}$, then
one must have $u \leq v$, $u \leq q_n$, and $v \leq p_n$.

Before describing examples which verify the inequalities above, we
note a few relevant facts.  First, if $f = \sum_j 2^j \chi_{E_j}$,
where the $E_j$ are pairwise disjoint measurable sets, and if $1
\leq p,u \leq \infty$, then
\[
    \Lp{f}{p,u} \sim
    (\sum_j2^{ju}|E_j|^{\frac{u}{p}})^{\frac{1}{u}},
\]
where the implicit constant depends on $p$ and $u$.  Second, if
$0<\vareps,r<1$, then $|N_{\vareps,r}| \sim
\vareps^{n-1}r^{n(n+1)/2}$ and $|B_{\vareps,r}| \sim
\vareps|N_{\vareps,r}|$.  Moreover, if $0<\vareps,r<1$, and $x \in
\reals^n$, then $T \chi_{N_{\vareps,r}(x)} \sim r$ on
$B_{\vareps,r}(x)$, so $\scriptt(N_{\vareps,r}(x),B_{\vareps,r}(x))
\sim r|B_{\vareps,r}| \sim
|N_{\vareps,r}|^{1/p_n}|B_{\vareps,r}|^{1/q_n'}$.

For the inequality $u \leq v$, we let $a = n+1$, and for
$j=1,2,\ldots$, we define $E_j = N_{2^{-aj}}(x_j)$, and $G_j =
B_{2^{-aj}}(x_j)$, where the $x_j$ are chosen so that the $E_j$, and also the $G_j$, are pairwise disjoint.  Then if
$f=\sum_1^M 2^{\frac{n-1}{p_n}aj}\chi_{E_j}$ and $g = \sum_{1}^M
2^{\frac{n}{q_n'}aj}\chi_{G_j}$, one has
\begin{align*}
    \Lp{f}{p_n,u} \sim M^{1/u} \qquad  \Lp{g}{q_n',v'}  \sim M^{1/v'}
\end{align*}
and $\langle Tf,g\rangle \gtrsim M$.  Thus for $T$ to map $L^{p_n,u}$ to $L^{q_n,v}$ boundedly, we must
have $M \leq M^{1/u+1/v'}$ for all positive integers $M$, i.e. $u
\leq v$.

We motivate our next two examples as follows.  Let a positive
integer $M$, positive constants $c,\eta$, and a set $\scriptj$ of
integers with $\#\scriptj = M$ be fixed.  Suppose that for each $j
\in \scriptj$ we have a pair $E_j,F_j$ of Borel sets so that
\begin{align*}
    \scriptt(E_j,F_j) \sim |E_j|^{1/p_n} |F_j|^{1/q_n'},\\
    2^{jp_n}|E_j| \sim c,\, \text{ and } \, |F_j| \sim \eta,
\end{align*}
where the implicit constants depend only on the dimension $n$ and
the exponents $p_n,q_n$.  Suppose further that the $E_j$ are
pairwise disjoint, as are the $F_j$.  Let $f=\sum_{\scriptj}
2^j\chi_{E_j}$ and $F = \bigcup_{\scriptj} F_j$.  Then
\begin{align*}
    \langle Tf, \chi_F \rangle &= \sum_{\scriptj} 2^j \scriptt(E_j,F_j) 
    \sim \sum_{\scriptj} 2^j|E_j|^{1/p_n}|F_j|^{1/q_n'} \sim M c^{1/p_n} \eta^{1/q_n'}.
\end{align*}
On the other hand, $\Lp{f}{p,u} \sim M^{1/u}c^{1/p_n}$ and
$|F|^{1/q_n'} \sim M^{1/q_n'}\eta^{1/q_n'}$.  Therefore, if such a
construction is possible for each positive integer $M$, boundedness
of $T:L^{p_n,u} \to L^{q_n,v}$ for any $v$ implies that $u \leq q_n$.  The
above construction could conceivably be used in a more general
context to produce counter-examples to Lorentz space bounds for any
operator with a rich enough family of quasi-extremals.  This general
construction is due to Christ (personal communication).

Now we construct a specific counter-example to demonstrate the
necessity of $u \leq q_n$.  Using our estimates on $|E_j| =
|N_{\vareps_j,r_j}(x_j)|$ and $|F_j| = |B_{\vareps_j,r_j}(x_j)|$, we
see that for $2^{jp_n}|E_j| \sim c$ and $|F_j| \sim \eta$, we must
have $r_j \sim \eta^{2/n(n+1)}\vareps_j^{-2/(n+1)}$, and $\vareps_j
\sim 2^{jp_n}\eta c^{-1}$.  If we let $\eta=2^{-Mn(n+1)/2}$,
$c=2^{Mp_n}\eta$, and $\scriptj=\{1,\ldots,M\}$, then $\vareps_j =
2^{jp_n-Mp_n}$ and $r_j = 2^{-j}$, which are both less than or equal
to 1 when $j \in \scriptj$.  Now choosing the sequence $x_j$ so that
the $E_j$, and likewise the $F_j$, are pairwise disjoint, we have
our counter-example.

The verification of the inequality $v \geq p_n$ is similar.  Now we
let $E=\bigcup_{\scriptj}E_j$ and $g=\sum_{\scriptj}2^j\chi_{F_j}$,
where $|E_j| \sim \eta$ and $2^{jq_n'}|F_j| \sim c$.  Again taking
$E_j = N_{\vareps_j,r_j}(x_j)$ and $F_j = B_{\vareps_j,r_j}(x_j)$,
where the $x_j$ will be chosen so the $E_j$ (and the $F_j$) are
pairwise disjoint, we compute $r_j \sim
\vareps_j^{-1/q_n}\eta^{2/n(n+1)}$ and $\vareps_j \sim c
\eta^{-1}2^{-jq_n'}$.  If we let $\scriptj = \{-1,\ldots,-M\}$,
$\eta = 2^{-Mq_n'(n-1)}$, and $c=2^{-Mq_n'n}$, then $\vareps_j =
2^{-(M+j)q_n'}$ and $r_j = 2^{jq_n'/q_n}$, which are both less than
or equal to 1 when $j \in \scriptj$.  The necessity of $v' \leq
p_n'$ follows by arguments similar to those two paragraphs above.

\section{A reduction to two multilinear bounds}
\label{sec:reduction}

In this section we state a theorem, essentially due to
Christ in \cite{QEx}, which will allow us to pass from a sort of
multilinear bound on characteristic functions of sets to the
strong-type inequality.  Let $S$ be
be a linear operator, mapping characteristic functions of Borel sets
to non-negative Borel functions.  If $E$ and $F$ are Borel
sets, then we define $\scripts (E,F):=\langle S\chi_E,\chi_F
\rangle$.

\begin{hyp1*}
    If $E_1$, $E_2$, and $F$ are Borel sets with positive, finite measures,
    if for $j=1,2$, $S
\chi_{E_j} \geq \alpha_j$ on $F$ and
    $\frac{\mathcal{S}(E_j,F)}{|E_j|} \geq \beta_j$,
    and if $\alpha_2 \geq \alpha_1$ and
$\beta_1\geq\beta_2$,
    then there exist real numbers $u_1$, $u_2$,
    $u_3$, and $u_4$, taken from a finite list depending on $S$ and satisfying
    $u_1+u_2 = \frac{r'}{r'-s'}$, $u_3+u_4 = \frac{r}{s-r}$,
and
    $\frac{u_2}{r} - \frac{u_4}{r'}-1>0$, such that
    \[
        \alpha_1^{u_1}\alpha_2^{u_2}\beta_1^{u_3}\beta_2^{u_4}
        \lesssim |E_2|,
    \]
    where the implicit constant depends on $S$ alone.
\end{hyp1*}
\begin{hyp2*}
    If $E$, $F_1$, and $F_2$ are Borel sets with positive finite measures,
    if for $j=1,2$, $S^*\chi_{F_j} \geq \beta_j$ on $E$ and
    $\frac{\mathcal{S}(E,F_j)}{|F_j|} \geq \alpha_j$, and if
    $\alpha_1\geq\alpha_2$ and
$\beta_2 \geq \beta_1$,
    then there exist real numbers $v_1$, $v_2$,
    $v_3$, and $v_4$ taken from a finite list which depends only on
    $S$ and satisfying
$v_1+v_2 = \frac{s'}{r'-s'}$, $v_3+v_4 = \frac{s}{s-r}$, and
    $\frac{v_4}{s'} - \frac{v_2}{s} -1>0$, such that
    \[
        \alpha_1^{v_1}\alpha_2^{v_2}\beta_1^{v_3}\beta_2^{v_4}
        \lesssim |F_2|,
    \]
    where the implicit constant depends only on $S$.
\end{hyp2*}

\begin{theorem}\label{lemma:wktost}
Let $S$ be a linear operator, mapping characteristic functions of
Borel sets to non-negative Borel measurable functions. Let $r$ and
$s$ be real numbers with $1 < r < s < \infty$, and $u$ and $v$ be
real numbers with $u<s$, $u<v$, and $r<v$. Suppose that Hypothesis 1
and Hypothesis 2 hold. Then the operator $S$ extends to a bounded
linear operator from $L^{r,u}(\reals^n)$ to $L^{s,v}(\reals^n)$.
\end{theorem}

As a partial motivation for the specific form of the hypotheses, we
initially observe that Hypothesis 2 is simply Hypothesis 1 for the
operator $S^*$ and exponents $(s',r')$ instead of $(r,s)$ (which is not to say that the hypotheses are equivalent).  Secondly,
under either hypothesis,  $S$ is of restricted weak type $(r,s)$, as
can be seen by letting $E_1=E_2$ in Hypothesis 1 or $F_1=F_2$ in
Hypothesis 2. Indeed, if $E$ and $F$ are Borel sets having positive
finite measures, we let $F_0 = \left\{x \in F:S\chi_E(x) \geq
\frac{\scripts(E,F)}{2|F|}\right\}$.  Letting $\alpha =
\frac{\scripts(E,F)}{2|F|}$ and $\beta = \frac{\scripts(E,F_0)}{|E|}
\sim \frac{\scripts(E,F)}{|E|}$, we then have $\alpha^{u_1+u_2}
\beta^{u_3+u_4} \lesssim |E|$.  Substituting in the values of
$\alpha$ and $\beta$ and using our identities for the $u_i$, we have
$\scripts(E,F) \leq C |E|^{\frac{1}{r}}|F|^{\frac{1}{s'}}$, where
$C$ depends only on $S$.

Our main use for the multilinear inequalities will be to show that
under certain assumptions on the various sets involved, namely
`quasi-extremality' of the pairs $(E_j,F)$ (see \cite{QEx}), disjointness of $E_1$ and $E_2$, and dissimilarity of $|E_1|$ and $|E_2|$, Hypothesis 1 implies that
$E_1$ and $E_2$ interact strongly (via $S$) with nearly
disjoint subsets of $F$.  This will allow us to treat $S$ as roughly diagonal when it is applied to functions of the form $\sum_j 2^j \chi_{E_j}$; see the appendix for details.

For the operator $T$ considered in this paper, $T$ and $T^*$ are essentially the same operator, so we give a little more explanation as to why we cannot expect to avoid verifying the hypotheses separately.  First, as noted above, under the hypotheses, strong interaction of $E_1$ and $E_2$ with the same set implies that $E_1$ and $E_2$ have comparable sizes.  But there are two natural ways of characterizing the interaction between $E_j$ and $F$ as strong for an operator $S=S^*:L^r \to L^s$:  either $\scripts(E_j,F) \sim \eps|E_j|^{1/r}|F|^{1/s'}$ (this is the situation in Hypothesis 1), or $\scripts(E_j,F) \sim \eps |F|^{1/r}|E_j|^{1/s'}$ (as in Hypothesis 2, with $E$ and $F$ exchanged).  {\it{A priori}}, there is no reason for these different types of strong interaction to have the same outcome.  Second, as will be seen in the appendix, Hypothesis 1 implies that $S$ is of weak-type $(r,s)$, while Hypothesis 2 implies that $S^*$ is of weak type $(s',r')$ (or $S$ is of restricted strong-type $(r,s)$).  These statements are not equivalent in general when $s \neq r'$.

The proof of Theorem \ref{lemma:wktost} essentially amounts to
changing exponents in \S 8 of \cite{QEx} and the addition of an extra hypothesis to handle the case when $r \neq s'$; we will give a complete proof in the appendix.

\section{The multilinear inequalities} \label{sec:multilinear}

In this section, we prove that Hypotheses 1 and 2 do hold for the
operator $T$ when $(r,s) = (p_n,q_n)$, where $T,p_n,q_n$ are as in
the introduction.  By Theorem \ref{lemma:wktost} and interpolation with the $L^1 \to L^1$ and $L^{\infty} \to L^{\infty}$ bounds, this will establish the main theorem.  Lemmas \ref{lemma:mlF} and Lemma \ref{lemma:mlE} verify Hypotheses 2 and 1, respectively.  We state the more complicated of the two lemmas first.

\begin{lemma} \label{lemma:mlF}
Assume that $E$, $F_1$, $F_2$ are Borel subsets of $\reals^n$ with
finite positive measures. Assume that $T^*\chi_{F_j}(x)
\geq \beta_j$ for $x \in E$ and that
$\frac{\mathcal{T}(E,F_j)}{|F_j|} \geq \alpha_j$.  Then
if $\alpha_2\leq\alpha_1$ and $\beta_2\geq\beta_1$,
\begin{align}\label{eq:mlF}
    |F_2| &\gtrsim
    \alpha_1^{r_1}\alpha_2^{r_2}\beta_1^{s_1}\beta_2^{s_2},
\end{align}
for some integers $r_j$ and $s_j$ (taken from a finite list which
depends on $n$), which satisfy
\begin{align}
    \label{eq:r1r2} \frac{n(n-1)}{2}&= r_1+r_2 ,\\
    \label{eq:s1s2}  n&=s_1+s_2, \text{\hspace{.5cm}and}\\
    \label{eq:ineqexp} 0&<\frac{s_2}{q_n'}-\frac{r_2}{q_n}-1.
\end{align}
\end{lemma}

\begin{lemma} \label{lemma:mlE}
Assume that $E_1$, $E_2$, $F$ are Borel subsets of $\mathbb{R}^n$
with finite positive measures.  Assume that
$T\chi_{E_j}(x)\geq \alpha_j$ for $x \in F$ and that
$\frac{\mathcal{T}(E_j,F)}{|E_j|}\geq \beta_j$.  Then if
$\alpha_2 \geq \alpha_1$, we have
\begin{align}\label{eq:mlE}
    |E_2| &\gtrsim \alpha_2^n
    \alpha_1^{\frac{n(n-1)}{2}}(\frac{\beta_1}{\alpha_1})^{n-1}.
\end{align}
\end{lemma}

We will comment on the differences between these lemmas at the end of this section.

\subsection{Proof of Lemma \ref{lemma:mlF}} \label{sec:pf_mlF}
First we prove Lemma \ref{lemma:mlF}. The main difficulty here is in
satisfying requirement (\ref{eq:ineqexp}), which is needed in the
proof of the strong-type inequality.

Let $\Phi_k:[-1,1]^k \to \mathbb{R}^n$ be defined by
\[\Phi_k(t) =h(t_1)-h(t_2)
+h(t_3)-\ldots +(-1)^{k+1}h(t_k).
\]
By Lemma 1 in \cite{CCC}, there exist a constant $c_n>0$ (which we
may assume is as small as needed), a point $x_0 \in E$, and Borel
sets $\Omega_k \subset[-1,1]^k$ for $1 \leq k \leq 2n-2$ such that
the following hold: $\Omega_{k+1}\subset \Omega_k \times [-1,1]$,
$|\Omega_1| = c_n \beta_1$, for each odd $k \leq 2n-3$ and $t \in
\Omega_k$,
\begin{itemize}
\item[]$x_0+\Phi_k(t) \in F_1$,
\item[]$|t_k-t_j| \geq c_n \beta_1$ if $j < k$,
\item[] and $|\{s \in [-1,1]:(t,s) \in \Omega_{k+1}\}| = c_n
\alpha_1$,
\end{itemize}
\noindent and for each even $k$ and $t \in \Omega_k$,
\begin{itemize}
\item[]$x_0 + \Phi_{k}(t) \in E$,
\item[]$|t_k-t_j|\geq c_n \alpha_1$ if $j<k$,
\item[] and if $k < 2n-2$, $|\{s \in [-1,1]:(t,s) \in
\Omega_{k+1}\}| = c_n \beta_1$.
\end{itemize}

Since $T^* \chi_{F_2}(x) \geq \beta_2$ on $E$, provided
$c_n$ is small enough ($<\frac{1}{2n}$ will do), there exists a
Borel set $\Omega_{2n-1}\subset \Omega_{2n-2} \times[-1,1]$ such
that if $t' \in \Omega_{2n-2}$, $|\{s \in [-1,1]:(t',s)
\in \Omega_{2n-1}\}| =c_n\beta_2$, and if $t \in \Omega_{2n-1}$,
then first, $x_0 + \Phi_{2n-1}(t) \in F_2$ and second,
$|t_{2n-1}-t_j|\geq c_n \beta_2$ whenever $j<2n-1$.

If $\beta_1 \gtrsim \alpha_1$, our lower bound is almost immediate,
and is essentially Lemma 2 in \cite{CCC} (there proved when
$F_1=F_2$). Fix $t^0 \in \Omega_{n-1}$ and let $\omega = \{t =
(t_1,t_2,\ldots,t_n) \in [-1,1]^n:(t^0,t) \in
\Omega_{2n-1}\}$. Then
\begin{align*}
|\omega| &\gtrsim
    \alpha_1^{\left\lfloor \frac{n}{2}\right\rfloor}
    \beta_1^{\left\lceil \frac{n}{2}\right\rceil}
    \frac{\beta_2}{\beta_1}.
\end{align*}
If we let $J(t)=|\det \frac{\partial
\Phi_{2n-1}(t^0,t)}{\partial t}|$, then $J(t)$
is just the absolute value of the Vandermonde determinant,
$J(t) = a_n \prod_{1 \leq i < j \leq n} |t_j-t_i|$, where
$a_n>0$.  By our lower bounds on $|t_i-t_j|$ when $i < j$ and $t \in
\Omega_j$,
\begin{align*}
 J(t) &\gtrsim
\begin{cases}
    \text{$(\frac{\beta_2}{\beta_1})^{n-1}
    \beta_1^{1+3+\cdots+n-1}\alpha_1^{2+4+\cdots+n-2}$ if $n$ is
    even}\\
    \text{$(\frac{\beta_2}{\beta_1})^{n-1}
    \beta_1^{2+4+\cdots+n-1}\alpha_1^{1+3+\cdots+n-2}$ if $n$ is
    odd}.
\end{cases}
\end{align*}

By an argument in \cite{CCC} using Bezout's theorem (see \cite{Shaf},
\cite{CCC}),
\begin{align*}
    |F_2| & \gtrsim \int_{\omega}|J(t)|dt 
    \gtrsim
    \begin{cases}
    \text{$(\frac{\beta_2}{\beta_1})^n
    (\frac{\beta_1}{\alpha_1})^{2+4+\cdots+n}
    \alpha_1^{\frac{n(n+1)}{2}}$ if $n$ is even}\\
    \text{$(\frac{\beta_2}{\beta_1})^n
    (\frac{\beta_1}{\alpha_1})^{1+3+\cdots+n}
    \alpha_1^{\frac{n(n+1)}{2}}$ if $n$ is odd}
    \end{cases}\\
    &\gtrsim (\frac{\beta_2}{\beta_1})^n
    (\frac{\beta_1}{\alpha_1})^n \alpha_1^{\frac{n(n+1)}{2}} =
    \beta_2^n \alpha_1^{\frac{n(n-1)}{2}},
\end{align*}
where the last inequality follows from our assumption that $\beta_1
\gtrsim \alpha_1$.  This is (\ref{eq:mlF}), with
$r_1=\frac{n(n-1)}{2}$, $r_2=s_1=0$, and $s_2=n$.  We
note that when $n=2$, we have $|F_2| \gtrsim \beta_2^2 \alpha_1$,
which implies Hypothesis 1, whether or not $\beta_1 \gtrsim
\alpha_1$.

In the case $\beta_1 < \alpha_1$, we extend the `band structure'
arguments of \cite{CCC}.  We will construct a partition of the integers
$\{1,2,\ldots,2n-1\}$ and use this partition
to pick out
$n$-dimensional subsets, or `slices', of $\Omega_{2n-1}$ such that
the Jacobian of the restriction of $\Phi_{2n-1}$ to these slices is
large.

Suppose that a partition of $\{1,2,\ldots,2n-1\}$ into subsets,
called `bands', is
given.  We designate each of the indices $1,2,\ldots,2n-1$ as free,
quasi-free, or bound as follows:
\begin{itemize}
    \item[] The least element of each band is free,
    \item[] If a band contains exactly two elements, the greater is
    quasi-free, and is quasi-bound to the lesser
    \item[] If a band contains three or more elements, the elements
    which are not least are designated as bound to the least element
    of that band.
\end{itemize}

Let $0<\vareps<1$ be fixed for now; it will be chosen (depending on
$n$ alone) to satisfy the hypotheses of a coming lemma.  We will
actually construct two partitions, the second a refinement of the
first.  In the first partition, 1 and all of the even indices will
be designated as free, and we will choose parameters
$0<c_{n,\varepsilon}<\delta'<\varepsilon\delta$ and a subset
$\omega$ of $\Omega_{2n-1}$ with $|\omega| \sim |\Omega_{2n-1}|$
such that for each $t \in \omega$,
\begin{itemize}
    \item[] $|t_i-t_j| \geq \delta \alpha_1$ unless $i$ and $j$
    belong to the same band
    \item[] If $i$ is quasi-bound to $j$, then $c_n\beta_1
    \leq |t_i-t_j| < \delta \alpha_1$
    \item[] If $i$ is bound to $j$, then $|t_i-t_j| <
    \delta' \alpha_1$.
\end{itemize}
We note that without the requirement (\ref{eq:ineqexp}), Lemma
\ref{lemma:mlF} could be proved by arguments in \cite{CCC} using
only this initial partition.  Inequality (\ref{eq:ineqexp}),
however, is an essential ingredient of our proof of Theorem
\ref{lemma:wktost} (used to prove for instance that the sets
(\ref{eq:gj}) are essentially disjoint).

Let $\scriptb$ be the band created in the first partition which
contains $2n-1$ and ignore, for now, all of the other bands and the
designation of the elements of $\scriptb$ as free, quasi-free, or
bound from the first partition. The second partition will subdivide
$\scriptb$. We will choose parameters
$c_{n,\varepsilon}<\rho'<\varepsilon \rho < \delta'$ and a subset
$\omega'$ of $\omega$ with $|\omega'| \sim |\omega|$ such that for
each $t \in \omega'$ and $i,j \in \scriptb$,
\begin{itemize}
    \item[] $|t_i-t_j| \geq \rho \gamma_2$ unless $i$ and $j$ belong
    to the same band
    \item[] If $i$ is quasi-bound to $j$, then
    $c_n \beta_1 \leq |t_i-t_j| < \rho \gamma_2$
    \item[] If $i$ is bound to $j$, then $|t_i-t_j|
    \leq \rho' \gamma_2$,
\end{itemize}
where $\gamma_2 = \max \{\alpha_2, \beta_2\}$.  

In \cite{CCC}, Christ developed an algorithm which, when $F_1=F_2$,
produces a band structure having the properties we want for the
first step of our partition.  For us, $F_1 \neq F_2$ in general, so
to achieve the first partition we must modify his argument.  The
second step, the refinement of the first partition, though a new
ingredient, will be achieved simply by a second application of the
algorithm--again with modifications in place. Though the needed
changes to the argument in \cite{CCC} are minor, for clarity we will
present the details of the algorithm in full.

Our algorithm will involve several iterative refinements of
certain sets.  To simplify exposition, if $\omega$ is a Borel set,
$\omega'$ and $\omega''$ will always denote Borel sets contained in
$\omega$ with $|\omega'|, |\omega''| \gtrsim |\omega|$, where the
implicit constant depends on $n$ alone (sometimes we will explicitly
specify this constant, sometimes not).

Initially set $\omega = \Omega_{2n-1}$.  Then there exist $\omega'
\subset \omega$ and a permutation $\sigma$ of $\{1,2,\ldots,2n-1\}$
such that $|\omega'| \geq \frac{1}{(2n-1)!}|\omega|$ and such that
$t\in\omega'$ and $i<j$ implies that
$t_{\sigma(i)}<t_{\sigma(j)}$. We henceforth
denote $\omega'$ by $\omega$.

Temporarily set $\delta = \frac{c_n}{2n}$ and $\delta' =
\frac{\varepsilon}{2}\delta$.

There exist $\omega' \subset \omega$, a positive integer $R$, and a
sequence of integers $1=L_1 < L_2 < \ldots < L_R \leq 2n-1$ such
that for each point $t \in \omega'$, $t_{\sigma(j)}-t_{\sigma(j-1)}
\geq \delta \alpha_1$ if and only if $j=L_i$ for some $1 < i \leq
R$.  To see this, note that it is possible to choose such an integer
$R$ and sequence for each $t \in \omega$ and that there are only
finitely many such sequences.  With no loss of generality, we may
assume $\omega'=\omega$.

Consider the following partition of $\{1,2,\ldots,2n-1\}$:
\begin{align*}
    \{\sigma(1)=\sigma(L_1),\sigma(L_1+1),&\ldots,\sigma(L_2-1)\}  \cup 
    \ldots \cup \{\sigma(L_R),
    \sigma(L_R+1),\ldots,\sigma(2n-1)\}.
\end{align*}
{\it A priori}, if $i$ and $j$ are in the same band,
$|t_i-t_j|<(2n-2)\delta \alpha_1 < c_n \alpha_1$. Since
$|t_j-t_i|\geq c_n\alpha_1$ when $j$ is even and $i<j$, each even
integer must be the least element of its band.  Therefore $R\geq n$
and no band has more than $n$ elements.  So if $i$ and $j$ are in
different bands, $|t_i-t_j| \geq \delta \alpha_1$, 
and if $i$ is quasi-bound to $j$, $c_n\beta_1 \leq |t_i-t_j| < \delta \alpha_1$
(for the lower bound, note that for $t \in \Omega_{2n-1}$ and $i
\neq j$, our assumptions $\alpha_1>\beta_1$ and $\beta_2>\beta_1$
imply that $|t_i-t_j|\geq c_n\beta_1$).  These two bounds are good,
but if $i$ is bound to $j$, we only know $|t_i-t_j| \leq (n-1)\delta \alpha_1$.

If for some $\omega' \subset \omega$ with $|\omega'|\geq
\frac{1}{2}|\omega|$, $t \in \omega'$ implies $|t_i-t_j|<\delta'
\alpha_1$ whenever $i$ is bound to $j$, we may assume $\omega =
\omega'$ and have the desired partition.  If there is no such
$\omega'$, then there is an $\omega'' \subset \omega$ and an index
$i_0$ bound to an index $j_0$ such that $|t_{i_0}-t_{j_0}| \geq
\delta' \alpha_1$ whenever $t \in \omega''$. In this case we start
over, setting $\omega = \omega''$, and selecting a new integer and
sequence as above, with $(\delta,\delta')$ replaced by
$(\frac{1}{n}\delta',\frac{\varepsilon}{2n}\delta')$
(with the latter pair now denoted by $(\delta,\delta')$).

Since our new parameters are less than our old, our old sequence of
indices will be a subsequence of our new one, and since
$|t_{i_0}-t_{j_0}| \geq n\delta \alpha_1$ ($n\delta$ being equal to
our old parameter $\delta'$), $i_0$ and $j_0$ must lie in different
bands.  We repeat this process until we have
$|t_i-t_j|<\delta'\alpha_1$ whenever $i$ is bound to $j$, on a
subset $\omega'$ of $\omega$ with $|\omega'| \geq \half |\omega|$.
Since there is at least one new band after the second iteration, we
have increased the number of free indices by at least one.  If there
are no bound indices, then we have satisfied the requirements for
the first partition; hence the process must terminate after at most
$n$ repetitions.

This completes the first partition.  We then partition
$\mathcal{B}$, as specified a few paragraphs above, in a similar
manner.

Our bands will be the subsets (other than $\scriptb$) from the two
partitions. Since 1 and all of the even indices are free, initially
we have at least $n$ free indices. If the total number of free and
quasi-free indices is exactly $n$, we let $\tilde{\omega} =
\Omega_{2n-1}$.

Suppose that there are more than $n$ indices which are free or
quasi-free.  Then we simply throw away the index 1, and redesignate
the indices $\{2,3,\ldots,2n-1\}$ as free, quasi-free, or bound
according to the rules above.  This increases the number of free and
quasi-free indices by one if exactly two indices were bound to 1,
does not change the number of free and quasi-free indices if three
or more indices were bound to 1, and decreases the total number by
one if 1 was free with no indices bound to it (and one or no indices
quasi-bound to it).  In particular, the total number of indices
designated as free or quasi-free never decreases by more than 1, so
this process leaves at least $n$ indices still designated as free or
quasi-free. We continue, successively throwing away indices
$1,2,\ldots,j$ and redesignating the remaining indices
$\{j+1,j+2,\ldots, 2n-1\}$, until we have a total of $n$ free or
quasi-free indices.  This process terminates after at most $n-1$
repetitions since the number of free and quasi-free indices cannot
be greater than the total number of indices.

Assume the above process terminates after $2n-k-1$ repetitions.  Henceforth,
when we refer to indices as free, quasi-free, or bound, we will be
referring only to those indices in $\{2n-k,\ldots,2n-1\}$, and to
their designation after the process above.

Because $|\omega| \sim |\Omega_{2n-1}|$ and by our upper bounds on
$|\{s \in [-1,1]:(t,s) \in \Omega_{j+1}\}|$ for $t \in \Omega_j$,
we may choose $t^0
\in \Omega_{2n-k-1}$ such that
\begin{align*}
    |\tilde{\omega}| & \gtrsim
    \alpha_1^{\left\lfloor \frac{k}{2}\right\rfloor}
    \beta_1^{\kovrtwo}
    \frac{\beta_2}{\beta_1},
\end{align*}
where $\tilde{\omega} = \{t=(t_{2n-k},t_{2n-k+1},\ldots,t_{2n-1}) \in [-1,1]^k:
(t^0,t) \in \omega\}$.

Now we explain how the slices are chosen.  Let
$2n-k=j_1<j_2<\ldots<j_n\leq 2n-1$ be the free or quasi-free
indices.  For $t \in \tilde{\omega}$, let
$\tau(t)=(\tau_1,\ldots,\tau_n)=(t_{j_1},\ldots,t_{j_n})$.  Let
$2n-k<i_1<i_2<\ldots<i_{k-n}\leq 2n-1$ be the bound indices, where
$i_{\ell}$ is bound to the index $j_{B(\ell)}$.  For
$t\in\tilde{\omega}$, let $s(t)=(s_{i_1},s_{i_2},\ldots,
s_{i_{k-n}}) = (t_{i_1}-t_{j_{B(1)}},t_{i_2}-t_{j_{B(2)}},\ldots,
t_{i_{k-n}}-t_{j_{{B(k-n)}}})$.
Then $t \mapsto (\tau(t),s(t))$ has an
inverse (which is linear) that we denote by $t(\tau,s)$.  Let
\[
    J(\tau,s)=|\det \frac{\partial\Phi_{2n-1}(t^0,t(\tau,s))}{\partial
    \tau}|.
\]
The proof of Lemma 3 of \cite{CCC} implies the following:

\begin{lemma} \label{lemma:eps}
Assume that we are given a partition of the indices
$\{2n-k,\ldots,2n-1\}$, where $k\leq 2n-1$, into bands such that exactly
$n$ indices are free or quasi-free.  Then there exists $\varepsilon
> 0$, depending only on $n$, such that if $t:\reals^k \to \reals^k$
is defined as above, and if $|t_i-t_j| < \vareps|t_j-t_{\ell}|$
whenever $j$ and $\ell$ are distinct free or quasi-free indices and
$i$ is bound to $j$, then
\begin{align} \label{eq:det}
    J(\tau,s) \gtrsim \prod_{1\leq i<j\leq n}|\tau_i-\tau_j|.
\end{align}
\end{lemma}

We choose our parameters $\delta,\delta'$ and $\rho,\rho'$ in the
band structure algorithm so $\delta'<\vareps \delta$,
$\rho<\delta'$, and $\rho'<\vareps \rho$, where $\vareps$ satisfies
the hypotheses of the lemma.  Thus if $j$ is free, $i$ is bound to $j$,
and $\ell$ is free or quasi-free and $\neq j$, then we have $|t_i-t_j|<\vareps|t_j-t_{\ell}|$ whenever $t \in \tilde{\omega}$.
Therefore (\ref{eq:det}) holds whenever $t(\tau,s) \in
\tilde{\omega}$.

Now we are finally ready to prove the lower bound on $|F_2|$ from Lemma \ref{lemma:mlF}.  We
consider first the case when the index $2n-1$ is free.

First suppose that $\mathcal{B} \cap \{2n-k,\ldots,2n-1\} =
\{2n-1\}$.  Then for each index $j$ with $2n-k \leq j<2n-1$, the
quantity $|t_{2n-1}-t_j|$ is bounded below by $\delta \alpha_1$
in addition to $c_n \beta_2$.  Hence by (\ref{eq:det}), for
$t(\tau,s) \in
\tilde{\omega}$,
\begin{align*}
    J(\tau,s) & \gtrsim \frac{\beta_2}{\alpha_1}\alpha_1^{n-1}
    \prod_{1\leq i < j \leq n-1}|\tau_j-\tau_i| 
    \gtrsim \frac{\beta_2}{\alpha_1}\alpha_1^{n-1}
    \alpha_1^{\frac{(n-1)(n-2)}{2}} (\frac{\beta_1}{\alpha_1})^M,
\end{align*}
where $M$ is the number of quasi-free indices.  By an application of
Bezout's theorem (as in \cite{CCC}), for each $s$,
\begin{align*}
    |F_2| & \gtrsim \int_{\{\tau:t(\tau,s) \in \tilde{\omega}\}}
    J(\tau,s) d\tau 
     \gtrsim \frac{\beta_2}{\alpha_1} \alpha_1^{\frac{n(n-1)}{2}}
    (\frac{\beta_1}{\alpha_1})^M |\{\tau:t(\tau,s) \in
    \tilde{\omega}\}|.
\end{align*}
For each $s \in \reals^{k-n}$ such that $t(\tau,s) \in
\tilde{\omega}$ for some $\tau \in \reals^n$, we must have
$|s_i|<\delta' \alpha_1$ for each $i$.  Therefore, integrating with
respect to $s$ on both sides of the inequality above,
\[
    \alpha_1^{k-n}|F_2| \gtrsim \frac{\beta_2}{\alpha_1}
    \alpha_1^{\frac{n(n-1)}{2}} (\frac{\beta_1}{\alpha_1})^M
    |\tilde{\omega}|.
\]
Thus
\begin{align*}
    |F_2| &\gtrsim \frac{\beta_2}{\alpha_1}\frac{\beta_2}{\beta_1}
    \alpha_1^{\frac{n(n+1)}{2}} (\frac{\beta_1}{\alpha_1})^M 
    (\frac{\beta_1}{\alpha_1})^{\kovrtwo}.
\end{align*}
If $k$ is even, the $\frac{k}{2}$ even indices and $2n-1$ are free.
Since the total number of free and quasi-free indices equals $n$,
$M+\frac{k}{2}+1 \leq n$.  If $k$ is odd, the $\frac{k-1}{2}$ even
indices, the index $2n-k$ (being the least index), and the index
$2n-1$ are free, so
$M+\frac{k+1}{2} + 1 \leq n$. In either case, our assumption that
$\beta_1 < \alpha_1$ implies
\[
    |F_2| \gtrsim \frac{\beta_2}{\alpha_1}\frac{\beta_2}{\beta_1}
    \alpha_1^{\frac{n(n+1)}{2}}(\frac{\beta_1}{\alpha_1})^{n-1}.
\]
Using $q_n=\frac{n(n+1)}{2(n-1)}$, one can
immediately check that the various equalities and inequality
(\ref{eq:ineqexp}) concerning the exponents
$r_1=\frac{n(n-1)}{2}$, $r_2=0$, $s_1=n-2$, and $s_2=2$
hold for $n \geq 3$.

We note in particular that if $\beta_2 \geq \alpha_1$, then
$|t_{2n-1}-t_j| \geq c_n\beta_2 \geq c_n \alpha_1 \geq \delta
\alpha_1$, so we are in the preceding case.  Henceforth, we will
assume that $\beta_2 < \alpha_1$, and hence that $\gamma_2 <
\alpha_1$.

Now to complete our analysis of the case when the index $2n-1$ is
free, we suppose that, in addition to $2n-1$, $\mathcal{B} \cap
\{2n-k,\ldots,2n-1\}$ contains at least one other free index.  We let
$M_1$ denote the number of quasi-free indices which are not
contained in $\mathcal{B}$, $M_2$ the number of quasi-free indices
which are contained in $\mathcal{B}$, and $N$ the number of free and
quasi-free indices in $\mathcal{B}$.

We then have that for $t(\tau,s) \in \tilde{\omega}$,
\[
    J(\tau,s) \gtrsim \alpha_1^{\frac{n(n-1)}{2}}(\frac{\gamma_2}
    {\alpha_1})^{\frac{N(N-1)}{2}-M_2}(\frac{\beta_1} {\alpha_1}
    )^{M_1+M_2}.
\]
Now if $\{\tau:t(\tau,s) \in \tilde{\omega}\} \neq \emptyset$,
$|s_i| \leq \delta' \alpha_1$ for each $i$ which is bound but not in
$\mathcal{B}$ and $|s_i| \leq \rho'\gamma_2$ for each bound index
$i$ in $\mathcal{B}$.  Therefore by Bezout's theorem and integration
in the possible values of $s$, which lie in a cube of size
$\sim \alpha_1^{k-n}(\frac{\gamma_2}{\alpha_1}) ^{R_2}$ ($R_2$ being the number of bound indices in $\mathcal{B}$),
\[
    |F_2| \gtrsim \alpha_1^{\frac{n(n+1)}{2}} (\frac{\gamma_2}
    {\alpha_1})^{\frac{N(N-1)}{2}-M_2-R_2}(\frac{\beta_1} {\alpha_1}
    )^{M_1+M_2} \frac{\beta_2}{\beta_1}
    (\frac{\beta_1}{\alpha_1})^{\kovrtwo}.
\]
Since
$2n-1$ is free, our arguments for the case when $\mathcal{B} =
\{2n-1\}$ imply that $M_1+M_2+ \kovrtwo < n$. Therefore, since $\gamma_2
< \alpha_1$,
\begin{align*}
    |F_2| &\gtrsim
    \alpha_1^{\frac{n(n+1)}{2}}
    (\frac{\beta_2}{\alpha_1})^{n-M_1-M_2-\kovrtwo}
    (\frac{\beta_1}{\alpha_1})^{M_1+M_2+ \kovrtwo}
    (\frac{\alpha_2}{\alpha_1})^
    {\max\{0,\frac{N(N-1)}{2}+M_1+\kovrtwo-n\}}
    \frac{\beta_2}{\beta_1}.
\end{align*}
The various requirements  on the exponents $r_1,r_2,s_1,s_2$ are
immediate except for (\ref{eq:ineqexp}).  If the max above is in
fact 0, the inequality follows from
\begin{align*}
    s_2=
    n-M_1-M_2- \kovrtwo+1
    \geq 2
\end{align*}
and $r_2=0$.

We note in particular that if $n=3$, by our assumption, the index 5
is free.  Therefore the index 3 must be free or quasi-free (because
0 or at least 2 indices may be bound), so we must have $k=3$, and 3
is actually free.  Hence $N = 2$ and $M_1=0$, so the max above is
zero.

If the max is not zero, we must show
\begin{align*}
    0&<
    n-\frac{M_2}{q_n'}-M_1-\frac{1}{q_n}-\frac{N(N-1)}{2q_n}-\kovrtwo.
\end{align*}
Since each quasi-free index in $\mathcal{B}$ is quasi-bound to a
unique free index in $\mathcal{B}$ and since $2n-1$ has no indices
quasi-bound to it, we have $M_2 \leq \frac{N-1}{2}$. If $k$ is even,
at most one of the $\frac{k}{2}$ even indices in
$\{2n-k,\ldots,2n-1\}$ is in $\mathcal{B}$.  Since the total number
of free and quasi-free indices is $n$, $M_1+\frac{k}{2} \leq n-N+1$.
A similar argument implies $M_1+\frac{k+1}{2} \leq n-N+1$ when $k$
is odd. Therefore the right hand side of the inequality above is
bounded from below by
\begin{align} \label{2inb}
    \frac{N-1}{2}-\frac{1}{q_n}-\frac{(N-1)^2}{2q_n}.
\end{align}

By the strict concavity of this term in $N$, it suffices to check
its positivity at the extreme values of $N$.  In general, we have $2
\leq N$, by assumption, and $N \leq n-\kovrtwo +1 \leq
\frac{n}{2}+1$ (since $k \geq n$, and at least $\kovrtwo +N -1$
indices are free or quasi-free).  When $n \geq 4$, one can check
that (\ref{2inb}) is positive at each of these points.

If we had $\gamma_2=\beta_2$, then the index $2n-1$ would be free,
since $|t_{2n-1}-t_j|\geq c_n \beta_2$ whenever $t \in
\Omega_{2n-1}$ and $j \neq 2n-1$.  Therefore, we may assume
henceforth that $\beta_2 < \gamma_2 = \alpha_2$.  We will not need
this fact for the case when $2n-1$ is quasi-free, but will need it
for the case when $2n-1$ is bound.

Suppose now that the index $2n-1$ is quasi-free.  Since
$|t_{2n-1}-t_j| \geq c_n \beta_2$ when $j<2n-1$ and $t \in
\tilde{\omega}$, we have for $t(\tau,s) \in \tilde{\omega}$
\[
J(\tau,s) \geq \alpha_1^{\frac{n(n-1)}{2}}
(\frac{\gamma_2}{\alpha_1})^{\frac{N(N-1)}{2}-M_2}
(\frac{\beta_1}{\alpha_1})^{M_1+M_2}
\frac{\beta_2}{\beta_1},
\]
where $M_1,M_2,N$ are as above.  Arguing as above,
\begin{align*}
&|F_2|  \gtrsim \alpha_1^{\frac{n(n+1)}{2}}
(\frac{\gamma_2}{\alpha_1})^{\frac{N(N-1)}{2}-M_2-R_2}
(\frac{\beta_1}{\alpha_1})^{M_1+M_2}
(\frac{\beta_2}{\beta_1})^2
(\frac{\beta_1}{\alpha_1})^{\kovrtwo}
\\
 &\gtrsim \alpha_1^{\frac{n(n+1)}{2}}
(\frac{\beta_2}{\alpha_1})^{n-M_1-M_2-\kovrtwo}
(\frac{\beta_1}{\alpha_1})^{M_1+M_2+ \kovrtwo}
 (\frac{\alpha_2}{\alpha_1})^
{\max\{0,\frac{N(N-1)}{2}+M_1+\kovrtwo-n\}}
(\frac{\beta_2}{\beta_1})^2,
\end{align*}
where now we know only that $n-M_1-M_2-\kovrtwo \geq 0$.  We check
that $0 < \frac{s_2}{q_n'}-\frac{r_2}{q_n}-1$.  This inequality holds if
the max above is zero, which by arguments similar to those above is
always the case when $n=3$.  Otherwise, $0 <
\frac{s_2}{q_n'}-\frac{r_2}{q_n}-1$ is equivalent to
\begin{align*}
    0 &<
    (n-M_1-M_2-\kovrtwo+2)(1-\frac{1}{q_n}) 
    -
    (\frac{N(N-1)}{2}+M_1+\kovrtwo-n)\frac{1}{q_n} -1.
\end{align*}
We simplify and use the bounds $M_2 \leq \frac{N}{2}$ and $n-N+1
\geq
    M_1+\kovrtwo$
to reduce the inequality (\ref{eq:ineqexp}) to showing
$\frac{N}{2}-\frac{2}{q_n} - \frac{N(N-2)}{2q_n}>0$. Whenever
$n \geq 4$, this can be proved via a concavity argument and the fact
that $2\leq N \leq \frac{n}{2}+1$ (the lower bound is because $2n-1$
and the index to which it is quasi-bound are in $\mathcal{B}$).

Finally we consider the case when the index $2n-1$ is bound.  Then
we have that $R_2$ (the number of bound indices in $\mathcal{B}$) is
at least 2.  Therefore
\begin{align*}
|F_2| &\gtrsim \alpha_1^{\frac{n(n+1)}{2}}
(\frac{\gamma_2}{\alpha_1})^{\frac{N(N-1)}{2}-M_2-2}
(\frac{\beta_1}{\alpha_1})^{M_1+M_2}
\frac{\beta_2}{\beta_1}
    (\frac{\beta_1}{\alpha_1})^{\kovrtwo} \\
& \geq \alpha_1^{\frac{n(n+1)}{2}} ( \frac{\beta_2}{\alpha_1}
)^{n-M_1-M_2- \kovrtwo } (
\frac{\alpha_2}{\alpha_1})^{\frac{N(N-1)}{2}+M_1+
\kovrtwo -n-2} ( \frac{\beta_1}{\alpha_1}
)^{M_1+M_2+  \kovrtwo }
\frac{\beta_2}{\beta_1},
\end{align*}
where the second inequality follows from $n-M_1-M_2
-\kovrtwo \geq 0$.  To establish (\ref{eq:mlF}) we must verify
that
\begin{align*}
0 &<(n-M_1-M_2-\kovrtwo+1)(1-\frac{1}{q_n})\\
&\hspace{1cm}-(\frac{N(N-1)}{2}+M_1+ \kovrtwo-n-2)(\frac{1}{q_n})-1.
\end{align*}
Since the index to which $2n-1$ is bound can have no indices
quasi-bound to it, $M_2 \leq \frac{N-1}{2}$.  Using the upper bound
for $M_1$ established above, the inequality will hold if
$\frac{N-1}{2}-\frac{(N-1)^2}{2q_n}+\frac{1}{q_n}>0$, which
holds for $N=1$ and $N=\frac{n}{2}+1$ and hence for all possible
values of $N$, for all $n \geq 3$.  This completes the proof of Lemma \ref{lemma:mlF}.

\subsection{Proof of Lemma \ref{lemma:mlE}.}

Now that the band structure argument has been described, the proof
of Lemma \ref{lemma:mlE} is much easier.  By the arguments above,
there exist a positive number $c_n$, a point $x_0 \in E_1$, and
measurable sets $\Omega_j \subset [-1,1]^j$ for $1 \leq j \leq 2n$
such that $|\Omega_1| = c_n \beta_1$ and $\Omega_{j+1} \subset
\Omega_j \times [-1,1]$ such that for each {\bf odd} $j \geq 1$ and
each $t \in \Omega_j$, $x_0+\Phi_j(t) \in F$, $|\{s \in [-1,1]:(t,s)
\in \Omega_{j+1}\}|=c_n\alpha_1$ if $j<2n-1$ and $=c_n\alpha_2$ if
$j=2n-1$, and $|t_j-t_i| \geq c_n \beta_1$ whenever $i <j$, and such
that for each {\bf even} $j<2n$ and each $t \in \Omega_j$,
$x_0+\Phi_j(t) \in E_1$, $|\{s \in [-1,1]:(t,s) \in \Omega_{j+1}\}|
= c_n \beta_1$, and $|t_j-t_i| \geq c_n \alpha_1$ whenever $i<j$,
and such that for each $t \in \Omega_{2n}$, $x_0+\Phi_{2n}(t)\in
E_2$ and $|t_{2n}-t_i| \geq c_n \alpha_2$ whenever $i<2n$.

We handle the case $\beta_1 \gtrsim \alpha_1$ in essentially the
same way that we did above.  Namely, let $t^0 \in \Omega_n$ and let
$\tilde{\omega} = \{t \in [-1,1]^n:(t^0,t) \in \Omega_{2n}\}$.  Then
$|\tilde{\omega}| \sim \alpha_1^{\left\lceil\frac{n}{2}\right\rceil}
\beta_1^{\left\lfloor\frac{n}{2}\right\rfloor}
\frac{\alpha_2}{\alpha_1}$.  Letting $J(t) = |\det \frac{\partial
\Phi_{2n}(t^0,t)}{\partial t}|$, for $t \in \tilde{\omega}$,
\[
J(t)\gtrsim
\begin{cases}
    \text{$\alpha_2^{n-1} \alpha_1^{1+3+\cdots+n-3}\beta_1^{2+4+\cdots
    n-2}$ if $n$ is even} \\
    \text{$\alpha_2^{n-1}
    \alpha_1^{2+4+\cdots+n-3}\beta_1^{1+3+\cdots+n-2}$ if $n$ is
    odd}.
    \end{cases}
\]
Hence
\begin{align*}
    |E_2| &\gtrsim \int_{\tilde{\omega}}J(t)dt 
    \gtrsim
    \begin{cases}
    \text{$\alpha_2^n \alpha_1^{n(n-1)/2} ( \frac{\beta_1}{\alpha_1}
    )^{1+3+\cdots+n-1}$, if $n$ is even}\\
    \text{$\alpha_2^n \alpha_1^{n(n-1)/2} (
    \frac{\beta_1}{\alpha_1} )^{2+4+\cdots+n-1}$ if $n$ is
    odd.}
    \end{cases}
\end{align*}
If $n=2$, $n=3$, or $\beta_1 \geq \alpha_1$, this implies the
inequality
\[
    |E_2| \gtrsim \alpha_2^n \alpha_1^{n(n-1)/2} (
    \frac{\beta_1}{\alpha_1})^{n-1}
\]
in Lemma \ref{lemma:mlE}.

Now supposing $\beta_1<\alpha_1$, we choose a subset $\omega \subset
\Omega_{2n}$, parameters $0<\delta'<\vareps\delta<\vareps c_n$,
where $\vareps$ satisfies the hypotheses of Lemma \ref{lemma:eps},
and a partition of the integers $\{1,\ldots,2n\}$ such that
$|\omega| \sim |\Omega_{2n}|$ and such that for each $t \in \omega$,
$|t_i-t_j| \gtrsim \delta \alpha_1$ unless $i$ and $j$ lie in the
same band, $c_n \beta_1 < |t_i-t_j|<\delta \alpha_1$ whenever $i$ is
quasi-bound to $j$, and $|t_i-t_j|<\delta' \alpha_1$ whenever $i$ is
bound to $j$.  Then by our assumption that $\alpha_2>\alpha_1$, the
index $2n$ is free.  Further, since $\delta<c_n$, 1 and every even
index other than $2n$ is also free.  Therefore we have at least
$n+1$ free or quasi-free indices after our initial partition.  We
proceed as above, dropping initial indices $1,2,\ldots, 2n-k$ and
redesignating the remaining indices until we have a total of $n$
free and quasi-free indices remaining in $\{2n-k+1,\ldots,2n\}$.  We
choose $t^0 \in \Omega_{2n-k}$ so that if we define $\tilde{\omega}
= \{t \in [-1,1]^k:(t^0,t) \in \omega\}$, we have $|\tilde{\omega}|
\gtrsim \alpha_1^{\kovrtwo}
\beta_1^{\left\lfloor\frac{k}{2}\right\rfloor}
\frac{\alpha_2}{\alpha_1}$.  Defining $J(\tau,s)$ as above, for
$t(\tau,s) \in \tilde{\omega}$,
\begin{align*}
    J(\tau,s) &\gtrsim \prod_{2n-k+1 \leq i<j\leq 2n}|t_i-t_j|
    \gtrsim (\frac{\alpha_2}{\alpha_1})^{n-1}
    \alpha_1^{\frac{n(n-1)}{2}}
    (\frac{\beta_1}{\alpha_1})^M,
\end{align*}
where $M$ is the number of quasi-free indices.  By Bezout's theorem,
for each $s \in \reals^{k-n}$
\begin{align*}
    |E_2| &\gtrsim \int_{\{\tau:t(\tau,s) \in \tilde{\omega}\}}J(\tau,s)
    d\tau 
    \gtrsim (\frac{\alpha_2}{\alpha_1})^{n-1}
    \alpha_1^{\frac{n(n-1)}{2}}
    (\frac{\beta_1}{\alpha_1})^M |\{\tau:t(\tau,s) \in
    \tilde{\omega}\}|.
\end{align*}
Integrating over the possible values of $s$, all of which lie in the
$k-n$ dimensional cube of diameter $2\delta'\alpha_1$,
\begin{align*}
    |E_2| & \gtrsim
    \alpha_1^{n-k}(\frac{\alpha_2}{\alpha_1})^{n-1}
    \alpha_1^{\frac{n(n-1)}{2}}
    (\frac{\beta_1}{\alpha_1})^M|\tilde{\omega}| 
    = (\frac{\alpha_2}{\alpha_1})^n
    \alpha_1^{\frac{n(n+1)}{2}}
    (\frac{\beta_1}{\alpha_1})^{M+\left\lfloor
    \frac{k}{2}\right\rfloor} \\
    & \gtrsim (\frac{\alpha_2}{\alpha_1})^n
    \alpha_1^{\frac{n(n+1)}{2}}
    (\frac{\beta_1}{\alpha_1})^{n-1},
\end{align*}
where the last line follows from the facts that $M+ \left\lfloor
\frac{k}{2}\right\rfloor \leq n-1$ (if $k$ is even, there are at
least $\frac{k}{2}+1$ free indices, and if $k$ is odd, there are at
least $\frac{k+1}{2}$ free indices) and $\beta_1 <\alpha_1$. This
proves Lemma \ref{lemma:mlE}.

\subsection{On the differences between Lemmas \ref{lemma:mlF} and \ref{lemma:mlE}}

Why, the reader may ask, do Lemmas \ref{lemma:mlF} and \ref{lemma:mlE} have different forms, and Lemma \ref{lemma:mlF} a more complicated proof?  Moreover, why is the quantity $\alpha_2$ involved in Lemma \ref{lemma:mlF} at all when $\beta_2$ is not involved in Lemma \ref{lemma:mlE} and $\alpha_2$ seems to play no role in the construction of $\Omega_{2n-1}$?

As mentioned in \S\ref{sec:reduction}, there are a few reasons to expect the hypotheses of Theorem \ref{lemma:wktost} to require separate verification.  Lemma \ref{lemma:mlF} is precisely Hypothesis 2, while Lemma \ref{lemma:mlE} is stronger than Hypothesis 1.  The statement of Lemma \ref{lemma:mlE} is stronger in large part because its proof is simpler.  We will give a few examples of some ``enemies" one encounters when trying simpler techniques than the two-stage band structure in the proof of Lemma \ref{lemma:mlF}.  

First, we explain why we cannot establish the necessary lower bound by stopping after the first stage of the band decomposition.  Say $\beta_1 \leq \beta_2 \ll \alpha_2 \leq \alpha_1$.  Assume that $2n-1$ is bound to another index.  Then the optimal lower bound on the Jacobian would be 
\[
|J(\tau,s)| \gtrsim \alpha_1^{n(n-1)/2} (\beta_1/\alpha_1)^M,
\]
for $t(\tau,s) \in \Omega$.  This implies the lower bound
\[
|F_2| \gtrsim \alpha_1^{n(n+1)/2} (\beta_1/\alpha_1)^{M+ \lceil k/2 \rceil} (\beta_2/\alpha_1),
\]
which is not strong enough to verify Hypothesis 2.

Another possibility might be to perform the band decomposition with $\beta_1$ as the quantity which determines whether an index is free or bound.  This is initially an attractive option in light of the fact that the inequality $\alpha_2 \geq \alpha_1$ is what makes the proof of Lemma \ref{lemma:mlE} so simple, whereas we have $\beta_2 \geq \beta_1$ in Lemma \ref{lemma:mlF}.  This we can also reject because in the proofs of Lemmas \ref{lemma:mlF} and \ref{lemma:mlE}, the band structure argument is only needed when $\beta_1 \leq \alpha_1$, in which case, the suggested decomposition would be trivial, and hence useless.

The author explored a few other possibilities for the decomposition, including choosing either $\alpha_2$ or $\beta_2$ instead of the parameter $\gamma_2=\max\{\alpha_2,\beta_2\}$, and found that none of these alternatives produced a strong enough lower bound.  This is not to say that no simpler alternative exists.

Finally, the issue of why $\alpha_2$ appears at all.  A glib but plausible answer would be that while $\alpha_2$ does not seem to play a role in the construction of $\Omega_{2n-1}$, it is a quantity which is intrinsic to the interaction between $E$ and $F_2$ and is hence relevant.  In fact, the actual identity of $\alpha_2$ plays no role in the proof, and we could just have well have taken any real number $0 < \rho \leq \alpha_1$ as ``our $\alpha_2$."  Because $r_2$ may be positive or negative, it is only by specializing to $\alpha_2$ that this alternative lemma implies Hypothesis 2.  That the statement of Lemma \ref{lemma:mlF} does not reflect the generality of the proof is a matter of aesthetics.  

\section{Proof of Theorem \ref{lemma:wktost}}
This proof is based on arguments due to Christ in \cite{QEx}.

We will henceforth assume that $r \leq u < v \leq s$.  We may do this with no loss of generality first because $t \geq \tilde{t}$ implies $\|f\|_{L^{p,t}} \leq \|f\|_{L^{p,\tilde{t}}}$ and second because increasing $u$ to $r$ or decreasing $v$ to $s$ if necessary adds no further restrictions to any of the exponents.  

We begin by proving that $S$ maps $L^{r,u} \to L^{s,\infty}$
boundedly. Along the way, we will prove an additional inequality
involving quasi-extremal pairs of sets.

By our assumptions on $S$, if $f$ and $g$ are functions with $f \geq
g$, then $Sf \geq Sg$.  Therefore it suffices to prove that
$\Lp{Sf}{s,\infty} \lesssim \Lp{f}{r,u}$ when $f$ is of the form
$\sum_j 2^j \chi_{E_j}$, where the $E_j$ are pairwise disjoint Borel
sets. Let $f = \sum_j 2^j \chi_{E_j}$ and let $F$ be a Borel set
having positive finite measure.  For each $\eta,\eps>0$ define
\begin{align*}
    \scriptj_0^F &= \{j \in \mathbb{Z}:
\scripts(E_j,F) = 0\},\\
    \scriptj_{\vareps}^F &= \{j \in \mathbb{Z}:
    \frac{\vareps}{2}|E_j|^{\frac{1}{r}}|F|^{\frac{1}{s'}} <
    \scripts(E_j,F) \leq
    \vareps |E_j|^{\frac{1}{r}}|F|^{\frac{1}{s'}}\}, \\
    \scriptj_{\vareps,\eta}^F &= \{j
\in \scriptj_{\vareps}^F:
    \frac{\eta}{2} < 2^{jr}|E_j|\leq\eta\}.
\end{align*}
Also for each $\eta$, $\vareps$, let
$\left\{\scriptj^F_{\vareps,\eta,i} \right\}_{i=1}^{\left\lceil A
\log \frac{1}{\vareps} \right\rceil}$ be a partition of
$\scriptj^F_{\vareps,\eta}$ into $A \left\lceil\log
\frac{1}{\vareps} \right\rceil$-separated subsets.  Here $A$ is a
large constant which will be determined later.

By the restricted weak-type bound, there exists a constant $C$ so
that
\begin{align}
\label{eq:unions}
    \mathbb{Z}\backslash \scriptj_0^F=  \bigcup_{0<\vareps \leq C}
    \scriptj_{\vareps}^F = \bigcup_{0<\vareps \leq C}
    \bigcup_{0<\eta \leq 1} \scriptj^F_{\vareps,\eta} = \bigcup_{0<\vareps\leq C}
    \bigcup_{0<\eta\leq 1}\bigcup_{i=1}^{\left \lceil A \log
    \frac{1}{\vareps} \right \rceil}\scriptj_{\vareps,\eta,i}^F,
\end{align}
where the union may be taken over dyadic values of $\eta$, $\vareps$.

Initially let $\vareps>0$ be fixed and assume that $\sum_{j \in
\scriptj_{\vareps}^{F}}2^{ju}|E_j|^{u/r} \leq 1$.  Let $\eta>0$ and
$1 \leq i \leq \left \lceil A\log\frac{1}{\vareps} \right\rceil$ be
fixed as well, and set $\scriptj = \scriptj^F_{\vareps,\eta,i}$.
Assume $\#\scriptj > 0$; otherwise $\scriptj = \emptyset$ and both
bounds below are trivial. We prove two bounds for $\sum_{j \in
\scriptj}2^j\scripts(E_j,F)$.

The first bound,
\begin{align*}
    \sum_{j \in \scriptj} 2^j \scripts(E_j,F) \sim
    \sum_{j \in \scriptj} \vareps 2^j |E_j|^{\frac{1}{r}}|F|^{\frac{1}{s'}} 
    \sim \vareps (\#\scriptj) \eta^{\frac{1}{r}}|F|^{\frac{1}{s'}} 
    \lesssim \vareps \eta^{\frac{1-u}{r}}|F|^{\frac{1}{s'}},
\end{align*}
follows from the restricted weak-type bound and our assumption that
$(\#\scriptj) \eta^{u/r} \sim \sum_{j \in \scriptj}2^{ju}|E_j|^{u/r}
\leq 1$.

For the second bound, given $j \in \mathbb{Z}$, let
\begin{align} \label{eq:gj}
    G_j = \left\{x \in F:S\chi_{E_j}(x) \geq
    \frac{\scripts(E_j,F)}{2|F|}\right\}.
\end{align}
By our assumption on $S$, $G_j$ is a Borel set. We will show that
Hypothesis 1 implies that $\sum_{j \in \scriptj}|G_j| \lesssim |F|$;
assume this for now.  Since $\scripts(E_j,F) \sim
\scripts(E_j,G_j)$,
\begin{align*}
    \sum_{j \in \scriptj} 2^j \scripts(E_j,F) &\lesssim
    \sum_{j \in \scriptj} 2^j |E_j|^{\frac{1}{r}}|G_j|^{\frac{1}{s'}}  \leq
    (\sum_{j \in \scriptj}2^{js}|E_j|^{\frac{s}{r}})^{\frac{1}{s}}
    (\sum_{j \in \scriptj}|G_j|)^{\frac{1}{s'}} \lesssim
    \eta^{(s-u)/rs}|F|^{\frac{1}{s'}},
\end{align*}
where the last two inequalities follow from H\"{o}lder's inequality.

We have just shown that
\begin{align}\label{ineq:jF,eps,eta,i}
    \sum_{j \in \scriptj^F_{\vareps,\eta,i}} 2^j \scripts(E_j,F)
    \lesssim \min\{\vareps \eta^{(1-u)/r}|F|^{\frac{1}{s'}},
    \eta^{(s-u)/(rs)}|F|^{\frac{1}{s'}}\}.
\end{align}
From (\ref{ineq:jF,eps,eta,i}) and (\ref{eq:unions}), if $1 \leq u<s$,
\begin{align*}
    \sum_{j \in \scriptj^F_{\vareps}} 2^j \scripts(E_j,F) & = \sum_{n=0}^{\infty} \sum_{i=1}^{\left\lceil
    A\log\frac{1}{\vareps}\right\rceil} \sum_{j \in \scriptj^F_{\vareps,2^{-n},i}}2^j\scripts(E_j,F) 
     \lesssim
        \vareps^a|F|^{\frac{1}{s'}},
\end{align*}
for some constant $a>0$ (since $s>u \geq 1$).  
If we drop the
requirement that $\sum_{j \in \scriptj_{\vareps}^F}2^{ju}|E_j|^{u/r}
\leq 1$, we then have
\begin{align} \label{eq:ineqqex}
    \sum_{j \in \scriptj^F_{\vareps}}2^j \scripts(E_j,F) \lesssim
    \vareps^a (\sum_{j \in \scriptj^F_{\vareps}}2^{ju}|E_j|^{u/r})^{\frac{1}{u}}
    |F|^{\frac{1}{s'}}.
\end{align}
Finally, by
(\ref{eq:unions}), we may sum over $\vareps = C2^{-m}$, for $0 \leq
m < \infty$, to obtain the bound (for any Borel measurable function
$f$)
\begin{align}\label{eq:ineqwk}
    \langle Sf,\chi_{F} \rangle \lesssim \Lp{f}{r,u}|F|^{\frac{1}{s'}}.
\end{align}

We will use (\ref{eq:ineqqex}) and (\ref{eq:ineqwk}) to prove that
$\Lp{Sf}{s,v} \lesssim \Lp{f}{r,u}$. We first note that we have
essentially the same inequalities for the operator $S^*$. Since
Hypothesis 2 is simply Hypothesis 1 with the operator $S$ replaced
by $S^*$ and $(r,s)$ replaced by $(s',r')$, if
$g = \sum_k 2^k \chi_{F_k}$, where the $F_k$ are pairwise disjoint
Borel measurable sets, if $E$ is Borel measurable, and if
$\scripts(E,F_k) \sim \vareps
|E|^{\frac{1}{r}}|F_k|^{\frac{1}{s'}}$ for each $k$, then whenever
$1 \leq v' < r'$,
\begin{align} \label{eq:ineqqex2}
    \sum_k 2^k \scripts(E,F_k) \lesssim \vareps^b |E|^{\frac{1}{r}}
    (\sum_{k}2^{kv'}|F_k|^{v'/s'})^{\frac{1}{v'}},
\end{align}
where $b$ is a positive constant. Moreover, if $g$ is a Borel
measurable function and $E$ a Borel measurable set, then
\begin{align} \label{eq:ineqwk2}
    \langle S\chi_E,g \rangle \lesssim |E|^{\frac{1}{r}}\Lp{g}{s',v'}.
\end{align}

Now for the strong-type bound, we assume $f = \sum_j2^j\chi_{E_j}$
and $g = \sum_k2^k\chi_{F_k}$, where the $E_j$, and likewise the
$F_k$, are pairwise disjoint Borel sets; we also assume that
$\Lp{f}{r,u} \sim 1$ and $\Lp{g}{s',v'} \sim 1$.  Let $0<\vareps
\leq 1$. Given $k$, we define $\scriptj_{\vareps}^{F_k}$ as above.
Given, in addition, $0<\eta \leq 1$, we define
\[
    \scriptk_{\eta} = \{k \in \mathbb{Z}:\frac{\eta}{2} < 2^{ks'}|F_k|
    \leq \eta\}.
\]
For $1 \leq i \leq \left\lceil A'\log(\frac{1}{\vareps})
\right\rceil$, we let $\scriptk_{\eta,i}^{\vareps}$ be an $A'\log
(\frac{1}{\vareps})$-separated subset of $\scriptk_{\eta}$,
where $A'$ will be chosen later.

We compute two bounds for $\sum_{k \in
\scriptk_{\eta,i}^{\vareps}}2^k \sum_{j \in
\scriptj_{\vareps}^{F_k}}2^j \scripts(E_j,F_k)$.

For the first bound, we use (\ref{eq:ineqqex}) to obtain
\begin{align*}
    \sum_{k \in \scriptk_{\eta,i}^{\vareps}}2^k
    \sum_{j \in \scriptj_{\vareps}^{F_k}}2^j \scripts(E_j,F_k)
    &\leq
    \sum_{k \in \scriptk_{\eta,i}^{\vareps}}2^k \vareps^a \Lp{f}{r,u}
    |F_k|^{\frac{1}{s'}} \\
    &\leq
    \vareps^a \Lp{f}{r,u} \eta^{(1-v')/s'} \Lp{g}{s',v'} \leq
    \vareps^a \eta^{(1-v')/s'},
\end{align*}
where we used the definition of $\scriptk_{\eta,i}^{\vareps}$ and
H\"{o}lder's inequality for the last line.

For the second bound, define
\begin{align} \label{eq:ejk}
    E_{j,k} = \left\{x \in E_j:S^*\chi_{F_k}(x) \geq
    \frac{\scripts(E_j,F_k)}{2|E_j|}\right\}.
\end{align}
Then $E_{j,k}$ is a Borel set.
We prove later that we may choose $A'$ so that for each
$j$, $\sum_{k \in \scriptk_{\eta,i}^{\vareps}:  j \in \scriptj_{\varepsilon}^{F_k}}|E_{j,k}| \lesssim
    |E_j|$;
assume this for now.  Then by (\ref{eq:ineqqex}) and H\"{o}lder's
inequality,
\begin{align*}
    \sum_{k \in \scriptk_{\eta,i}^{\vareps}} 2^k
    \sum_{j \in \scriptj_{\vareps}^{F_k}} 2^j \scripts(E_j,F_k)
        & \lesssim \sum_{k \in \scriptk_{\eta,i}^{\vareps}} 2^k
    (\sum_{j \in \scriptj_{\vareps}^{F_k}}
    2^{ju}|E_{j,k}|^{u/r}
    )^{\frac{1}{u}}|F_k|^{\frac{1}{s'}} \\
    & \leq (\sum_{k \in \scriptk_{\eta,i}^{\vareps}}
    \sum_{j \in \scriptj_{\vareps}^{F_k}} 2^{ju} |E_{j,k}|^{u/r})^{\frac{1}{u}}
    (\sum_{k \in \scriptk_{\eta,i}^{\vareps}}2^{ku'}
    |F_k|^{\frac{u'}{s'}} )^{\frac{1}{u'}} \\
    & \lesssim (\sum_{j\in\mathbb{Z}} 2^{ju}|E_j|^{\frac{u}{r}})^{1/u} \eta^{(u'-v')/(s'u')}
    (\sum_{k \in
    \scriptk_{\eta,i}^{\vareps}}2^{kv'}|F_k|^{v'/s'}
    )^{\frac{1}{u'}}\\ &\leq \eta^{(u'-v')/(s'u')}.
\end{align*}
Here the third inequality follows from Minkowski's inequality, our assumption that $u\geq r$, and the choice of $A'$ mentioned earlier in this paragraph.

Now, letting $A''=\left\lceil A'\log2 \right\rceil$, when $u<v$ we
have
\begin{align*}
    \langle Sf,g \rangle &= \sum_{j,k} 2^j2^k\scripts(E_j,F_k) 
    = \sum_{m=0}^{\infty}\sum_{n=0}^{\infty}
    \sum_{k \in \scriptk_{2^{-m}}}2^k \sum_{j \in
    \scriptj^{F_k}_{C2^{-n}}}2^j\scripts(E_j,F_j) \\
    &= \sum_{m=0}^{\infty}\sum_{n=0}^{\infty}
    \sum_{i=1}^{A''n} \sum_{k \in
\scriptk_{2^{-m},i}^{2^{-n}}}
    2^k \sum_{j \in
    \scriptj_{2^{-n}}^{F_k}}2^j \scripts(E_j,F_k) \\
    & \lesssim \sum_{m=0}^{\infty}\sum_{n=0}^{\infty}
    \sum_{i=1}^{A''n} \min\{2^{-m(u'-v')/(s'u')},2^{-na}
    2^{m(v'-1)/s'}\} \\
    &\lesssim 1,
\end{align*}
since $u<v$ implies $u'>v'$.

To complete the proof of the lemma, it remains to show that
if $G_j$ is defined as in (\ref{eq:gj}), then we may choose $A$
so that $\sum_{\scriptj_{\vareps,\eta,i}^{F}}|G_k| \lesssim |F|$,
and that if $E_{j,k}$ is defined as in (\ref{eq:ejk}),
it is possible to choose $A'$ so that $\sum_{k \in \scriptk_{\eta,i}^{\vareps}:
    j \in \scriptj_{\varepsilon}^{F_k}}|E_{j,k}| \lesssim  |E_j|$.
The two situations are essentially symmetric, so it suffices to
prove the former.

Let $\scriptj = \scriptj_{\vareps,\eta,i}^{F}$.  By H\"older's
inequality, we have that
\begin{align*}
    (\sum_{j \in \scriptj}|G_j|)^2 &=
    (\int_F \sum_{j \in\scriptj}\chi_{G_j})^2 \leq
    |F|\int(\sum_{j \in\scriptj}\chi_{G_j})^2 =
    |F|(\sum_{j \in \scriptj}|G_j|+\sum_{j \neq k \in \scriptj}
    |G_j \cap G_k|).
\end{align*}
Thus, either $\sum_{j \in \scriptj}|G_j| \lesssim |F|$ or
$(\sum_{j \in \scriptj}|G_j|)^2 \lesssim |F|\sum_{j \neq
k \in \scriptj}|G_j \cap G_k|$.  The former is the inequality we
want, so assume the latter occurs.  From the restricted weak type
bound on $S$ and our definitions of $\scriptj$ and $G_j$,
\[
    |E_j|^{\frac{1}{r}}|G_j|^{\frac{1}{s'}} \gtrsim \scripts(E_j,G_j) \sim
    \scripts(E_j,F) \gtrsim \vareps|E_j|^{\frac{1}{r}}|F|^{\frac{1}{s'}},
\]
for each $j \in \scriptj$.  Hence
\begin{align*}
    (\#\scriptj\vareps^{s'}|F|)^2 &\lesssim
    (\sum_{j \in \scriptj}|G_j|)^2 \lesssim
    |F|\sum_{j \neq k \in \scriptj}|G_j \cap G_k| \lesssim |F|(\#\scriptj)^2 \max_{j\neq k \in \scriptj}
    |G_j \cap G_k|,
\end{align*}
i.e. there exist distinct indices $j,k \in \scriptj$ such
that $|G_j \cap G_k| \gtrsim \vareps^{2s'}|F|$.

Assume $j>k$ and let $G = G_j \cap G_k$.
Since $G \subset G_j$, for $x \in G_j$,
\[
    S \chi_{E_j}(x) \gtrsim \frac{\scripts(E_j,F)}{|F|} \gtrsim
    \vareps|E_j|^{\frac{1}{r}}|F|^{\frac{-1}{s}} \gtrsim
    \vareps^{\frac{s-3}{s-1}}|E_j|^{\frac{1}{r}}|G|^{\frac{-1}{s}}=:\alpha_j.
\]
Similarly, $S \chi_{E_k}(x) \gtrsim
\vareps^{\frac{s-3}{s-1}}
|E_k|^{\frac{1}{r}}|G|^{\frac{-1}{s}}=:\alpha_k$ on $G$.  We also
have that
\begin{align*}
    \frac{\scripts(E_j,G)}{|E_j|} \gtrsim \vareps^{2s'+1}|E_j|^{\frac{-1}{r'}}
    |G|^{\frac{1}{s'}}&=:\beta_j \\
    \frac{\scripts(E_k,G)}{|E_k|} \gtrsim \vareps^{2s'+1}|E_k|^{\frac{-1}{r'}}
    |G|^{\frac{1}{s'}}&=:\beta_k.
\end{align*}
Now we use Hypothesis 1.  Since $j>k$, and $j,k \in \scriptj$,
$|E_j| < |E_k|$ (by the definition of $\scriptj^F_{\vareps,\eta}$),
so $\alpha_j<\alpha_k$ and $\beta_j > \beta_k$.  Therefore
$\alpha_j^{u_1}\alpha_k^{u_2}\beta_j^{u_3}\beta_k^{u_4}
\lesssim |E_k|$.
By our assumptions on the exponents $u_i$, this implies
\begin{align*}
    |E_k| &\gtrsim \vareps^{B_0}|E_j|^{\frac{u_1}{r}-\frac{u_3}{r'}}
    |E_k|^{\frac{u_2}{r}-\frac{u_4}{r'}}
    |G|^{\frac{u_3+u_4}{s'}-\frac{u_1+u_2}{s}} \\
        &= \vareps^{B_0}|E_j|^{1-\frac{u_2}{r}+\frac{u_4}{r'}}|E_k|^{\frac{u_2}{r}
    -\frac{u_4}{r'}},
\end{align*}
for a positive constant $B_0$ (independent of the $u_i$).
Now since $\frac{u_2}{r}-\frac{u_4}{r'}-1>0$, this implies that
$\vareps^{B_{u}}|E_k| \lesssim |E_j|$, where $B_{u}$ is
positive and depends on the $u_i$.  Since the $u_i$ are taken from a finite
list, we let $B$ be the maximum of the $B_{u}$, and let $A =
C\frac{B}{r \log 2}$ (C will depend on the implicit constant in the previous
sentence).  On the other hand, we are assuming that $\scriptj$ is $A \log
\frac{1}{\vareps}$-separated.  Since $|E_j| \sim 2^{-jr}\eta$ and
$|E_k| \sim 2^{-kr}\eta$, we have a contradiction.  Therefore we must have
$\sum_{j \in \scriptj}|G_j| \lesssim |F|$.


%

\end{document}